\documentclass[12pt]{amsart}

\usepackage{amsmath,amssymb,verbatim,graphicx}
\usepackage[mathscr]{eucal}
\usepackage{enumerate}

\newtheorem{theorem}{Theorem}
\newtheorem{conjecture}[theorem]{Conjecture}

\theoremstyle{definition}
\newtheorem{definition}{Definition}
\newtheorem{example}{Example}

\begin{document}

\title{A counterexample to Matsumoto's conjecture regarding 
absolute length vs. relative length in Finsler manifolds}

\author{Jeanne N. Clelland} 

\address{Dept. of Mathematics, 395 UCB, University of Colorado, Boulder, CO 80309-0395\\ Jeanne.Clelland@colorado.edu}

\thanks{The author was partially supported by a Collaboration Grant for Mathematicians from The Simons Foundation.}


\maketitle
\begin{abstract}
Matsumoto conjectured that for any Finsler manifold $(M, F)$ for which the 
restriction of the fundamental tensor to the indicatrix of $F$ is positive 
definite, the absolute length $F(X)$ of any tangent vector $X \in T_xM$ is 
the global minimum for the relative length $|X|_y$ as $y$ varies along the 
indicatrix $I_x \subset T_xM$ of $F$.  In this note, we disprove 
this conjecture by presenting a counterexample.
\end{abstract}

\section{Preliminaries}

Let $M$ be a smooth manifold of dimension $n$.  Let 
$x = (x^1, \ldots, x^n)$ be a local coordinate system on $M$, and let 
$(x, y) = (x^1, \ldots, x^n, y^1, \ldots, y^n)$ 
be the corresponding canonical local coordinate system  on the tangent bundle $TM$.
A {\em Finsler metric} on $M$ is a function
\[ F: TM \to [0, \infty) \]
with the following properties:
\vspace{-0.1in}
\begin{enumerate}
\item{Regularity: $F$ is $C^{\infty}$ on the slit tangent bundle $TM
\setminus \{0\}$.}
\item{Positive homogeneity: $F(x, \lambda y) = \lambda F(x, y)$ for
all $\lambda > 0$.}
\item{Strong convexity: The fundamental tensor $g = [ g_{ij} ]$, 
defined by
\[ g_{ij}(x,y) = \frac{1}{2}  \frac{\partial^2 (F^2)}{\partial y^i\,
\partial y^j} , \]
is positive definite at every point of $TM \setminus \{0\}$.}\label{strongcon}
\end{enumerate}
(For details, see \cite{BCS00}.)  In other words, a Finsler metric on
a manifold $M$ is a smoothly varying Minkowski norm on each tangent
space $T_xM$.  A Finsler metric is Riemannian if and only if the 
fundamental tensor $g$ is independent of the $y$-coordinates, in which case we have
\[ F^2(x,y) = g_{ij}(x)y^i y^j. \]
A {\em Finsler manifold} is a pair $(M, F)$ where $M$ is a 
smooth manifold and $F$ is a Finsler metric on $M$.

\begin{definition}
Let $(M, F)$ be a Finsler manifold.
\begin{itemize}
\item For any tangent vector $X \in T_xM$, the {\em length} of $X$, also 
known as the {\em absolute length} of $X$, is defined to be the value of $F(x, X)$.
\item For any $x \in M$, the {\em indicatrix} at $x$ is the subset $I_x \subset T_xM$ defined by
\[ I_x = \{y \in T_xM \mid F(x, y) = 1\}. \]
The indicatrix is the Finsler analog of the unit sphere in each tangent space.
\end{itemize}
\end{definition}

In \cite{Matsumoto79}, Matsumoto introduced the following notion of {\em 
relative length}:

\begin{definition}
For any tangent vector $X \in T_xM$ and any vector $y \in I_x$, the 
{\em relative length} $|X|_y$ of $X$ with respect to $y$ is defined by
\[ |X|_y = \sqrt{g_{ij}(x,y) X^i X^j}. \]
\end{definition}
It follows from the positive homogeneity property of Finsler metrics that 
for any $X \in I_x$, $|X|_X = F(X) = 1$.  More generally, for any nonzero 
$X \in T_xM$, if we set $\bar{X} = \frac{X}{F(X)}$, then $|X|_{\bar{X}} = F(X)$.

Based on this definition, Matsumoto defined the {\em relative energy} 
function of a given tangent vector $X \in T_xM$ to be the function 
$E_X:I_x \to \mathbb{R}$ defined by
\[ E_X(y) = \frac{1}{2}|X|_y^2 . \]
Matsumoto explored the behavior of the critical points of the relative 
energy function and made the following conjecture:

\begin{conjecture}\label{Matsumoto-conjecture}
For any nonzero $X \in T_xM$, the relative energy function 
$E_X:I_x \to \mathbb{R}$ achieves its global minimum at the point 
$y = \frac{X}{F(X)}$; equivalently, 
$F(X) \leq |X|_y$ for all $y \in I_x$.
\end{conjecture}

This conjecture was stated somewhat informally, with the caveat 
``provided that [the metric] satisfy the usual desirable assumptions."  It 
appears from the remainder of the paper \cite{Matsumoto79} that this 
comment refers to the strong convexity condition for Finsler metrics, and 
that Matsumoto also considered more general metrics for which this
condition does not necessarily hold everywhere.

\section{A counterexample}

In \cite{Matsumoto79}, Matsumoto proved Conjecture \ref{Matsumoto-conjecture} 
for {\em cubic} Finsler metrics, i.e., metrics of the form
\[ F(x,y) = \sqrt[3]{a_{ijk}(x) y^i y^j y^k}, \]
under the additional assumption that the induced metric of the indicatrix 
is positive definite.  (We note that this positivity condition is an 
immediate consequence of the strong convexity condition.)  However, the 
indicatrix $I_x$ of a 
metric of this form---or of an analogous $m$th-root metric for any odd 
$m$---in any tangent space $T_xM$ is never convex, or even compact, and so 
this positivity condition can never hold on all of $I_x$.  This also means that, strictly speaking, such a metric can never be Finsler.

\begin{figure}[h]
\begin{center}
\includegraphics[width=2in]{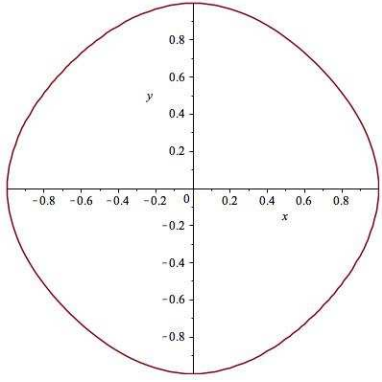}
\end{center}
\caption{Indicatrix for the metric $F(y) = \sqrt[4]{y_1^4 + 3 y_1^2 y_2^2 + y_2^4}$}
\label{indicatrix-fig}
\end{figure}

In this section, we present a counterexample to Conjecture \ref{Matsumoto-conjecture}. 
Let $M$ be a 2-dimensional Minkowski space with local coordinates 
$y = (y_1, y_2)$ (we use subscripts rather than superscripts in order to 
distinguish more easily between indices and exponents), and consider the 
$4$th-root metric
\begin{equation}\label{my-metric}
F(y) = \sqrt[4]{y_1^4 + 3 y_1^2 y_2^2 + y_2^4}
\end{equation}
on $M$.  The indicatrix of this metric is shown in Figure \ref{indicatrix-fig}.

The fundamental tensor of $F$ is readily computed to be
\begin{equation}\label{my-gij}
g = \frac{1}{2(y_1^4 + 3 y_1^2 y_2^2 + y_2^4)^{3/2}} 
\begin{bmatrix} 
2 y_1^6 + 9 y_1^4 y_2^2 + 6 y_1^2 y_2^4 + 3 y_2^6 & 5 y_1^3 y_2^3 \\[0.1in]
5 y_1^3 y_2^3 & 3 y_1^6 + 6 y_1^4 y_2^2 + 9 y_1^2 y_2^4 + 2 y_2^6
\end{bmatrix}.
\end{equation}
We have
\[ \text{trace}(g) = \frac{5\,(y_1^2 + y_2^2)^3}{2\,(y_1^4 + 3 y_1^2 y_2^2 + y_2^4)^{3/2}}, \qquad 
\det(g) = \frac{3\,(2 y_1^4 + y_1^2 y_2^2 + 2 y_2^4)}{4\,(y_1^4 + 3 y_1^2 y_2^2 + y_2^4)} , \]
and since both of these expressions are strictly positive away from the 
origin, it follows that $g$ is positive definite at every point of 
$M \setminus \{0\} $. Therefore, its restriction to the indicatrix $I$ of 
$F$ is positive definite as well.

Now let $X = (X_1, X_2) \in M$, and consider the relative energy function
\begin{equation}\label{my-energy}
\begin{aligned}
& E_X(y)  = \frac{1}{2} \left( g_{11}(y) X_1^2 + 2 g_{12}(y) X_1 X_2 + g_{22}(y) X_2^2 \right) \\[0.1in]
& = \frac{(2 X_1^2 + 3 X_2^2) y_1^6 + (9 X_1^2 + 6 X_2^2) y_1^4 y_2^2 
+ 10 X_1 X_2 y_1^3 y_2^3 + (6 X_1^2 + 9 X_2^2) y_1^2 y_2^4 
+ (3 X_1^2 + 2 X_2^2) y_2^6}
{4\, (y_1^4 + 3 y_1^2 y_2^2 + y_2^4)^{3/2}}.
\end{aligned}
\end{equation}
Note that the denominator of $E_X$ is constant on the indicatrix $I$. 
Therefore, in order to find critical points of $E_X$ on $I$, it suffices 
to find critical points of the numerator
\begin{equation}\label{my-energy-numerator}
\begin{aligned}
& \tilde{E}_X(y) = \\[0.1in]
& (2 X_1^2 + 3 X_2^2) y_1^6 + (9 X_1^2 + 6 X_2^2) y_1^4 y_2^2 
+ 10 X_1 X_2 y_1^3 y_2^3 + (6 X_1^2 + 9 X_2^2) y_1^2 y_2^4 
+ (3 X_1^2 + 2 X_2^2) y_2^6,
\end{aligned}
\end{equation}
subject to the constraint equation 
\[ F^4(y) = y_1^4 + 3 y_1^2 y_2^2 + y_2^4 = 1  \]
that defines the indicatrix.  To accomplish this, we use a Lagrange 
multiplier approach: We compute
\begin{equation} \label{my-gradients}
\begin{gathered} \nabla \tilde{E}_X = \begin{bmatrix}
y_1\left( (12 X_1^2 + 18 X_2^2)(y_1^4 + y_2^4) 
+ (36 X_1^2 + 24 X_2^2) y_1^2 y_2^2 + 30 X_1 X_2 y_1 y_2^3 \right) \\[0.1in]
y_2\left( (18 X_1^2 + 12 X_2^2)(y_1^4 + y_2^4) 
+ 30 X_1 X_2 y_1^3 y_2 + (24 X_1^2 + 36 X_2^2) y_1^2 y_2^2 \right)
\end{bmatrix}, \\[0.1in]
\nabla (F^4) = \begin{bmatrix}
y_1 (4 y_1^2 + 6 y_2^2) \\[0.1in]
y_2 ( 6 y_1^2 + 4 y_2^2)
\end{bmatrix}. 
\end{gathered}
\end{equation}
Critical points occur wherever the vectors $\nabla \tilde{E}_X$ and 
$\nabla(F^4)$ are linearly dependent.  A straightforward computation shows 
that this occurs precisely when
\begin{equation}\label{my-crit-point-condition}
y_1 y_2 (y_1 - y_2) (y_1 + y_2) (y_1^2 + y_2^2) (X_1 y_2 - X_2 y_1)^2 = 0.
\end{equation}
Therefore, the critical points $\bar{y} = (\bar{y}_1, \bar{y}_2) \in I$ 
are those points where any of the following hold:
\vspace{-0.1in}
\begin{itemize}
\item $\bar{y}$ is a scalar multiple of $X$;
\item $\bar{y}_1 = 0$, $\bar{y}_2 \neq 0$;
\item $\bar{y}_1 \neq 0$, $\bar{y}_2 = 0$;
\item $\bar{y}_1 = \pm \bar{y}_2$.
\end{itemize}
Evaluating the relative energy function \eqref{my-energy} at these points 
yields the critical values
\begin{equation}\label{my-crit-values}
\begin{aligned}
E_X(X)& = \sqrt{X_1^4 + 3 X_1^2 X_2^2 + X_2^4} = \frac{1}{2}F^2(X), \\
E_X(\bar{y}_1, 0) & = \frac{1}{2} X_1^2 + \frac{3}{4} X_2^2, \\
E_X(0, \bar{y}_2) & = \frac{3}{4} X_1^2 + \frac{1}{2} X_2^2, \\
E_X(\bar{y}_1, \bar{y}_1) & = \frac{1}{2\sqrt{5}}(2 X_1^2 + X_1 X_2 + X_2^2), \\
E_X(\bar{y}_1, -\bar{y}_1) & = \frac{1}{2\sqrt{5}}(2 X_1^2 - X_1 X_2 + X_2^2).
\end{aligned}
\end{equation}
(We note that, since the relative energy function \eqref{my-energy} is 
homogeneous of degree zero in $y$, it is not necessary to scale the 
critical points so that $F(\bar{y}) = 1$ in order to compute the critical values.)

\begin{example}
Let $X = (1,0)$.  Then from \eqref{my-crit-values}, we have
\[ E_X(X) = E_X(\bar{y}_1, 0) = \frac{1}{2}, \qquad E_X(0,\bar{y}_2) = \frac{3}{4}, \qquad E_X(\bar{y}_1, \pm \bar{y}_1) = \frac{1}{\sqrt{5}}, \]
so $E_X$ has neither a global maximum nor a global minimum at $X$. 

The graph of $E_X$ as a function of the Euclidean angle parameter 
$\theta = \tan^{-1}(\frac{y_2}{y_1})$ along $I$ is shown in Figure 
\ref{example1-fig}. The scalar multiples of $X$ in $I$ correspond to 
$\theta \in \{0, \pi, 2\pi\}$ and are indicated on the graph.

\begin{figure}[h]
\begin{center}
\includegraphics[width=2in]{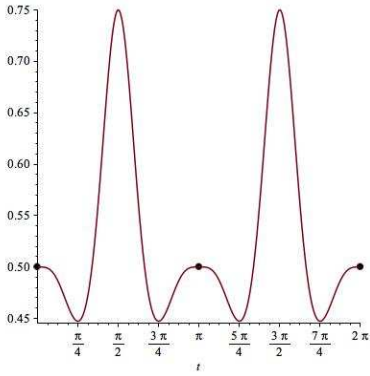}
\end{center}
\caption{Graph of $E_X$ for $X = (1,0)$}
\label{example1-fig}
\end{figure}

\end{example}

\begin{example}
For an example where $X$ is not one of the ``universal" critical points, 
let $X = (1,3)$.  Then from \eqref{my-crit-values}, we have
\begin{gather*}
E_X(X)  = \frac{\sqrt{109}}{2} \approx 5.22,  \qquad 
E_X(\bar{y}_1,0)  = \frac{29}{4} = 7.25, \qquad  
E_X(0, \bar{y}_2)  = \frac{21}{4} = 5.25, \\
E_X(\bar{y}_1, \bar{y}_1)  = \frac{23}{2\sqrt{5}} \approx 5.14, \qquad 
E_X(\bar{y}_1, \bar{y}_1) = \frac{1}{\sqrt{5}} \approx 3.80,
\end{gather*}
so again, $E_X$ has neither a global maximum nor a global minimum at $X$.

The graph of $E_X$ as a function of the Euclidean angle parameter 
$\theta = \tan^{-1}(\frac{y_2}{y_1})$ along $I$ is shown in Figure 
\ref{example2-fig}. The scalar multiples of $X$ in $I$ correspond to 
$\theta \in \{\tan^{-1}(3),\, \tan^{-1}(3) + \pi\}$ and are indicated on 
the graph.  In 
this case, the critical points at scalar multiples of $X$ are not even 
local extrema of $E_X$, but rather are inflection points.

\begin{figure}[h]
\begin{center}
\includegraphics[width=2in]{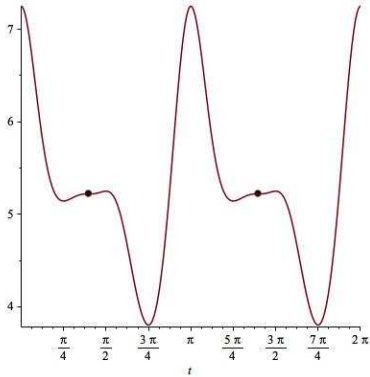}
\end{center}
\caption{Graph of $E_X$ for $X = (1,3)$}
\label{example2-fig}
\end{figure}

\end{example}

These computations clearly demonstrate that Conjecture 
\ref{Matsumoto-conjecture} does not hold for the metric 
\eqref{my-metric}, and hence is false in general.


\begin{thebibliography}{99}

\bibitem{BCS00} D. Bao, S.-S. Chern, and Z. Shen, {\em An
Introduction to Riemann-Finsler Geometry}, Graduate Texts in
Mathematics 200, Springer-Verlag, New York (2000).

\bibitem{Matsumoto79} M. Matsumoto, ``The length and the relative length of tangent vectors of Finsler spaces," {\em Rep.\ Math.\ Phys.\/} {\bf 15}  (1979) 375--386.


\end{thebibliography}
\end{document}